\documentclass[11pt]{amsart}

\usepackage{amsmath,amssymb,amsthm,amsfonts,amscd,flafter,epsf, epsfig,graphicx,verbatim,pinlabel,mathrsfs}
\usepackage[all]{xy}
\usepackage{epsf} 
\usepackage{epstopdf}

\title[Support norm]{A note on the support norm of a contact structure}
\author{John A. Baldwin}
\address{
    Department of Mathematics \\
    Princeton University}
\email{baldwinj@math.princeton.edu}
\urladdr{http://math.princeton.edu/\char126 baldwinj}

\author{John B. Etnyre}
\address{
    School of Mathematics\\
    Georgia Institute of Technology}
\email{etnyre@math.gatech.edu}
\urladdr{http://math.gatech.edu/\char126 etnyre}
\thanks{The first author was supported by an NSF Postdoctoral Fellowship and NSF Grant DMS-0635607. The second author was partially supported by NSF Grant DMS-0804820.}

\newtheorem{thm}{Theorem}
\newtheorem{lem}[thm]{Lemma}

\newtheorem{prop}[thm]{Proposition}
\newtheorem{quest}[thm]{Question}
\newtheorem{rem}[thm]{Remark}

\def\bn{\operatorname{bn}}
\def\sg{\operatorname{sg}}
\def\sn{\operatorname{sn}}
\def\rot{\operatorname{rot}}
\def\bfn{\mathbf{n}}
\def\bft{\mathbf{t}}

\begin{document}
\begin{abstract} 
In this note we observe that the no two of the three invariants defined for contact structures in \cite{et3} -- that is, the support genus, binding number and support norm -- determine the third. 
\end{abstract}

\maketitle

In \cite{et3}, the second author and B. Ozbagci define three invariants of contact structures on closed, oriented 3-manifolds in terms of supporting open book decompositions. These invariants are the support genus, binding number and support norm. There are obvious relationships between these invariants, but \cite{et3} leaves open the question of whether any two of them determine the third. We show in this note that this is not the case.

Recall that an open book decomposition $(L,\pi)$ of a 3--manifold $M$ consists of an oriented link $L$ in $M$ and a fibration $\pi: (M-L)\to S^1$ of the complement of $L$ whose fibers are Seifert surfaces for $L.$ The fibers $\pi^{-1}(\theta)$ of $\pi$ are called \emph{pages} of the open book and $L$ is called the \emph{binding}. It is often convenient to record an open book decomposition $(L,\pi)$ by a pair $(\Sigma,\phi)$, where $\Sigma$ is a compact surface which is homeomorphic to the closure of a page of $(L,\pi)$, and $\phi:\Sigma\rightarrow\Sigma$ is the monodromy of the fibration $\pi$. A contact structure $\xi$ on $M$ is said to be supported by the open book decomposition $(L,\pi)$ if $\xi$ is the kernel of a 1-form $\alpha$ which evaluates positively on tangent vectors to $L$ that agree with the orientation of $L$, and for which $d\alpha$ restricts to a positive volume form on each page of $(L,\pi)$.

With this in mind, we may describe the three invariants defined in \cite{et3}. The \emph{support genus} of a contact structure $\xi$ on $M$ is defined to be $$\sg(\xi)=\min\{g(\pi^{-1}(\theta))\,|\, (L, \pi) \text{ supports } \xi \},$$ where $\theta$ is any point in $S^1$ and $g(\pi^{-1}(\theta))$ is the genus of the page $\pi^{-1}(\theta)$. 
The \emph{binding number} of $\xi$ is defined to be $$\bn(\xi)=\min \{|L| \, | \, (L,\pi) \text{ supports } \xi \text{ and }  \sg(\xi)=g(\pi^{-1}(\theta))\},$$ 
where $|L|$ denotes the number of components of $L$ (or, equivalently, the number of boundary components of any page of $(L,\pi)$). And the \emph{support norm} of $\xi$ is defined to be $$ \sn(\xi)=\min\{-\chi(\pi^{-1}(\theta))\,|\, (L, \pi) \text{ supports } \xi \},$$
where $\chi(\pi^{-1}(\theta))$ denotes the Euler characteristic of any page $\pi^{-1}(\theta).$ It is a simple observation that $\sn(\xi)\geq -1$, with equality if and only if $\xi$ is the standard tight contact structure on $S^3.$

Since, for any surface $\Sigma$, we have the equality $$-\chi(\Sigma)=2g(\Sigma)+|\partial \Sigma |-2,$$ it is  immediately clear that $$ \sn(\xi)\leq 2\sg(\xi)+\bn(\xi)- 2.$$ Moreover, if the support norm of $\xi$ is achieved by an open book whose pages have genus $g>\sg(\xi)$ and whose binding has $m$ components, then $\sn(\xi) = 2g+m-2$, which is at least $2\sg(\xi)+1$. The following lemma from \cite{et3} summarizes these bounds.

\begin{lem}\label{bound}
For any contact structure $\xi$ on a closed, oriented 3--manifold, 
$$\min\{ 2\sg(\xi) + \bn(\xi) -2, \,2\sg(\xi) +1\} \,\leq \,\sn(\xi)\,\leq\, 2\sg(\xi) + \bn (\xi) -2.$$ 
\end{lem}

Thus, for contact structures with $\bn(\xi)\leq 3$, it follows that that $\sn(\xi)=2\sg(\xi) + \bn (\xi) -2.$ Yet, the results in \cite{et3} do not resolve whether the upper bound on the support norm in Lemma \ref{bound} can ever be a strict inequality. Our main result is that this bound can indeed be a strict inequality; that is, the support genus and binding number do not, in general, wholly determine the support norm.

\begin{figure}[!htbp]
\labellist 
\hair 2pt 
\small
\pinlabel $x$ at 65 148
\pinlabel $y$ at 101 85
\pinlabel $\delta$ at 122 30

\endlabellist 
\begin{center}
\includegraphics[height = 4.3cm]{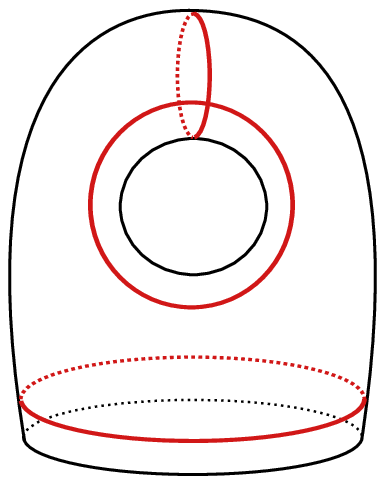}
\caption{\quad The surface $\Sigma$ and the curves $x$, $y$ and $\delta$.}
\label{fig:sfc}
\end{center}
\end{figure}

For the rest of this note, $\Sigma$ will denote the genus one surface with one boundary component. Let $\phi_{\bfn,m}$ be the diffeomorphism of $\Sigma$ given by $$\phi_{\bfn,m} = D_\delta^m \cdot D_{x}D_{y}^{-n_1}\cdots D_xD_y^{-n_k},$$ where $x$, $y$ and $\delta$ are the curves pictured in Figure \ref{fig:sfc}, and $\bfn = (n_1,\dots,n_k)$ is a $k$-tuple of non-negative integers for which some $n_i \neq 0.$ Let $\xi_{\bfn,m}$ denote the contact structure supported by the open book $(\Sigma,\phi_{\bfn,m})$, and let $M_{\bfn,m}$ denote the 3-manifold with this open book decomposition.

\begin{thm}\label{main}
For $m\leq 0$, the contact structure $\xi_{\bfn,m}$ satisfies $$\sg(\xi_{\bfn,m})=0.$$ For any fixed tuple $\bfn$, there is a finite subset $E_{\bfn}$ of the integers such that $$\bn(\xi_{\bfn,m})>3 \text{  and  } \sn(\xi_{\bfn,m}) = 1$$ for all $m\leq 0$ which are not in $E_{\bfn}$. In particular, $$\sn(\xi_{\bfn,m})< 2\sg(\xi_{\bfn,m})+\bn(\xi_{\bfn,m})-2$$ for all $m\leq 0$ which are not in $E_{\bfn}$.
\end{thm}

In contrast, the support genus $\sg(\xi_{\bfn,m})=1$ when $m>0,$ \cite{bald5}. Therefore, for $m>0$, $\bn(\xi_{\bfn,m}) = 1$ and $\sn(\xi_{\bfn,m}) = 1$, and, hence, $$\sn(\xi_{\bfn,m})=2\sg(\xi_{\bfn,m})+\bn(\xi_{\bfn,m})-2.$$ That is, the upper bound in Lemma \ref{bound} is achieved for $\xi=\xi_{\bfn,m}$ when $m>0.$ 

\begin{proof}[Proof of Theorem \ref{main}]
One can easily see that, for $m\leq 0$, the open book $(\Sigma,\phi_{\bfn,m})$ is not \emph{right-veering}; therefore, for such $m$, the contact structure $\xi_{\bfn,m}$ is overtwisted \cite{hkm1} and $\xi_{\bfn,m}$ is supported by an open book with planar pages \cite{et}.

Observe that if a contact structure $\xi$ on $M$ is supported by an open book with planar pages and the binding number of $\xi$ is three or less, then $M$ must be a Seifert fibered space. More precisely, if the binding number of $\xi$ is two then $M$ is a lens space, and if the binding number is three then $M$ is a small Seifert fibered space. One can see this by drawing a surgery picture corresponding to the open book supporting $\xi$ which realizes the binding number.

It is well known that the diffeomorphism $\phi_{\bfn,0}=D_{x}D_{y}^{-n_1}\cdots D_xD_y^{-n_k}$ is pseudo-Anosov (for instance, $\phi_{(1),0}$ is the monodromy of the figure eight knot in $S^3$). Therefore, the binding of the open book given by $(\Sigma, \phi_{\bfn,0})$ is a hyperbolic knot, and the manifold $M_{\bfn,m}$ is obtained from $M_{\bfn,0}$ via $-\frac 1m$ surgery on this knot. Thurston's Dehn Surgery Theorem then implies that there is some finite subset $E_{\bfn}$ of the integers for which $M_{\bfn,m}$ is hyperbolic for all $m$ not in $E_{\bfn}$ \cite{th1}. In particular, $M_{\bfn,m}$ is not a Seifert fibered space, save, perhaps, for some of the $m$ in the exceptional set $E_{\bfn}$. Hence, the binding number of $\xi_{\bfn,m}$ must be greater than three for all $m\leq 0$ which are not in $E_{\bfn}$. 

We are left to check that the support norm of $\xi_{\bfn,m}$ is one when $m\leq 0$ and $M_{\bfn,m}$ is hyperbolic. If the support norm were not one, then it would be zero (the support norm must be non-negative since $\xi_{\bfn,m}$ is not the tight contact structure on $S^3$). But the only surface with boundary which has Euler characteristic zero is the annulus, and the only 3-manifolds with open book decompositions whose pages are annuli are lens spaces. \end{proof}

It is natural to ask if the difference between $\sn(\xi)$ and $2\sg(\xi)+\bn(\xi)-2$ can be arbitrarily large. While we cannot answer this question we do note the following.
\begin{thm}\label{boundondifference}
For a fixed $\bfn$ the difference between $\sn(\xi_{\bfn, m})$ and $2\sg(\xi_{\bfn, m}) + \bn (\xi_{\bfn, m}) -2$ is bounded independent of $m<0.$
\end{thm}

Before we prove this theorem we estimate  the binding numbers $\bn(\xi_{\bfn,m})$ in some special cases. 

\begin{prop}
\label{bn}
The binding number of $\xi_{(1),-1}$ satisfies $3\leq\bn(\xi_{(1),-1})\leq 9.$ For each $m<-1$, the binding number of $\xi_{(1),m}$ satisfies $4\leq \bn(\xi_{(1),m})\leq 9.$
\end{prop}

The manifold $M_{(1),-1}$ is the Brieskorn sphere $\Sigma(2,3,7)$. Since $\Sigma(2,3,7)$ is not a lens space, it does not admit an open book decomposition with planar pages and two or fewer binding components. Therefore, $\bn(\xi_{(1),-1})\geq 3$. It is well-known that the only exceptional surgeries on the figure eight are integral surgeries \cite{th1}. Therefore, $E_{(1)} = \{-1,0,1\}$. So, from Theorem \ref{main}, we know that $\bn(\xi_{(1),m})>3$ for all $m<-1$. To prove Proposition \ref{bn}, we construct an open book decomposition of $M_{(1),m}$ with planar pages and nine binding components and we show that it supports $\xi_{(1),m}$ for $m<0$. 

Recall that overtwisted contact structures on a 3-manifold $M$ are isotopic if and only if they are homotopic as 2-plane fields. Moreover, the homotopy type of a 2-plane field $\xi$ is uniquely determined by its induced $Spin^c$ structure $\bft_{\xi}$ and its \emph{3-dimensional invariant} $d_3(\xi)$. Therefore, in order to show that the open book decomposition we construct actually supports $\xi_{(1),m}$ (for $m\leq 0$), we need only prove that the contact structure it supports is overtwisted and has the same 3-dimensional invariant as $\xi_{(1),m}$ (their $Spin^c$ structures automatically agree since $H_1(M_{(1),m};\mathbb{Z})=0$). Below, we describe how to compute these invariants from supporting open book decompositions. For more details, see the exposition in \cite{et3}.

%When $c_1(\bft_{\xi})$ is torsion in $H^2(M;\mathbb{Z})$, $d_3(\xi)$ is an element of $\mathbb{Q}$ and may be computed as follows \cite{gompf}. Let $(X,J)$ be an almost-complex 4-manifold whose boundary is $M$ and for which $\xi$ is homotopic as a 2-plane field to the complex tangencies along $M=\partial X$. Then, \begin{equation*}\label{eqn:d3}d_3(\xi) = \frac 14(c_1^2(X,J) - 2\chi(X)-3\sigma(X)).\end{equation*} 

Suppose that $\phi$ is a product of Dehn twists around homologically non-trivial curves $\gamma_1,\dots, \gamma_k$ in some genus $g$ surface $S$ with $n$ boundary components. The open book $(S,id)$ supports the unique tight contact structure on $\#^{2g+n-1}(S^1\times S^2)$, and the $\gamma_i$ may be thought of as Legendrian curves in this contact manifold. The contact manifold $(M,\xi)$ supported by the open book $(S,\phi)$ bounds an achiral Lefschetz fibration $X$, which is constructed from $\sharp^{2g+n-1}(S^1\times D^3)$ by attaching 2-handles along these Legendrian curves. Each 2-handle is attached with contact framing $\pm 1$ depending on whether the corresponding Dehn twist in $\phi$ is left- or right-handed, respectively. %If there are any left-handed Dehn twists in $\phi$, then $X$ does not admit an almost-complex structure. 
As long as $c_1(\bft_{\xi})$ is torsion in $H^2(M;\mathbb{Z})$, $d_3(\xi)$ is an element of $\mathbb{Q}$ and may be computed according to the formula, \begin{equation}\label{d32}
d_3(\xi)=\frac 14 (c^2(X)-2\chi(X)-3\sigma(X)) + q.
\end{equation}
Here, $q$ is the number of left-handed Dehn twists in the factorization $\phi.$ The number $c^2(X)$ is the square of the class $c(X) \in H^2(X;\mathbb{Z})$ which is Poincar{\'e} dual to $$\sum_{i=1}^k \rot(\gamma_i)C_i \in H_2(X,M;\mathbb{Z}),$$ where $C_i$ is the cocore of the 2-handle attached along $\gamma_i$, and $\rot(\gamma_i)$ is the rotation number of $\gamma_i$. The class $c(X)$ restricts to $c_1(\bft_{\xi})$ in $H^2(M;\mathbb{Z})$. Since we have assumed that $c_1(\bft_{\xi})$ is torsion, some multiple $k\cdot c(X)$ is sent to zero by the map $i^*:H^2(X;\mathbb{Z})\rightarrow H^2(M;\mathbb{Z})$, and, hence, comes from a class $c_r(X)$ in $H^2(X,M;\mathbb{Z})$, which can be squared. So, by $c^2(X)$, we mean $\frac 1 {k^2} c_r^2(X)$. 

%For a more detailed exposition of these constructions and formulae, see \cite{et3}; there, the authors also describe an easy way to compute the rotation numbers $r(\gamma_i)$.

%Before we describe the planar open book decomposition for $M_{(1),m}$ alluded to above, we compute the 3-dimensional invariants of the contact structures $\xi_{(1),m}$.

\begin{lem}
\label{d3hf}
For $m<0$, the 3-dimensional invariant $d_3(\xi_{(1),m}) = 1/2$.
\end{lem}

\begin{proof}[Proof of Lemma \ref{d3hf}]
Let $M_{(1),m}$ denote the result of $-\frac 1m$ surgery on the figure eight knot. It is therefore a rational homology 3-sphere, and the 3-dimensional invariant $d_3(\xi_{(1),m})$ is a well-defined element of $\mathbb{Q}$. Observe that the Dehn twist $D_{\delta}$ is isotopic to the composition $(D_xD_y)^6$. As described above, the contact manifold supported by the open book $(\Sigma,\phi_{(1),m})$ bounds an achiral Lefschetz fibration $X$, constructed from $\sharp^2(S^1\times D^3)$ by attaching $12|m|+2$ 2-handles corresponding to the Dehn twists in the factorization $$\phi_{(1),m} = (D_{x}D_y)^{6m} \cdot D_xD_y^{-1}.$$ From the discussion in \cite[Section 6.1]{et3}, it follows that $\rot(x) = \rot(y) = 0$; hence, $c(X)=0$. Moreover, $\chi(X) =  12|m|+1$ and $q=12|m|+1$. Therefore, the formula in (\ref{d32}) gives $$d_3(\xi_{(1),m})=\frac{12|m|+1}2 -\frac{3\sigma(X)} 4.$$ 

The achiral Lefschetz fibration associated to the monodromy $(D_xD_y)^{2m}$ gives a well-known Milnor fiber with the reverse orientation. Its signature is $8|m|.$ One may easily check (via Kirby calculus or gluing formulas for the signature or computations of the degree of related Heegaard-Floer contact invariants) that $\sigma(X)=8|m|.$ Thus $d_3(\xi_{(1),m)})=1/2.$
\end{proof}

\begin{figure}[!htbp]
\labellist 
\hair 2pt 
\small

\endlabellist 
\begin{center}
\includegraphics[height = 2.4cm]{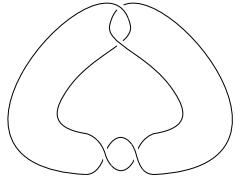}
\caption{\quad The figure eight, drawn here as a twist knot.}
\label{fig:fig82}
\end{center}
\end{figure}

\begin{figure}[!htbp]
\labellist 
\hair 2pt 
\small
\pinlabel $S$ at 20 90
\pinlabel $P$ at 910 90
\pinlabel $r$ at 1300 200
\tiny
\pinlabel $\gamma_1$ at 910 205
\pinlabel $\gamma_3$ at 1225 430
\pinlabel $\gamma_2$ at 1040 269
\pinlabel $\gamma_4$ at 1410 430

\endlabellist 
\begin{center}
\includegraphics[height = 4.3cm]{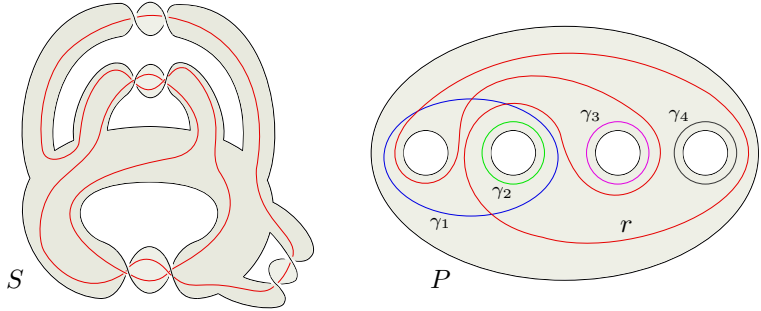}
\caption{\quad The figure eight knot, embedded on a planar surface with five boundary components.}
\label{fig:fig8}
\end{center}
\end{figure}

Proposition \ref{bn} follows if we can find a planar open book with nine binding components which supports an overtwisted contact structure on $M_{(1),m}$ with $d_3 = 1/2$. The figure eight knot $K$ is pictured in Figure \ref{fig:fig82}. $K$ can be embedded as a homologically non-trivial curve on the surface $S$ obtained by plumbing together two positive Hopf bands and two negative Hopf bands, as shown on the left in Figure \ref{fig:fig8}. 

Topologically, $S$ is an embedded copy of the planar surface $P$ with five boundary components shown on the right in Figure \ref{fig:fig8}. Moreover, $S$ is a page of the open book decomposition of $M_{(1),0}\cong S^3$ given by $(P,\phi),$ where $\phi$ is the product of right-handed Dehn twists around the curves $\gamma_1$ and $\gamma_3$ and left-handed Dehn twists around the curves $\gamma_2$ and $\gamma_4$. The knot $K$ is the image, under this embedding, of the curve $r\subset P$. 

Since the Seifert framing of $K$ agrees with the framing induced by $S$, $(P,D_r^{m}\cdot\phi)$ is an open book decomposition for $-\frac 1m$ surgery on $K$. Let $\xi_m$ denote the contact structure on $M_{(1),m}$ which is supported by this open book. It is easy to check that the open book $(P,D_r^{m}\cdot\phi)$ is not right-veering for $m\leq 0.$ (This can be seen by taking, for example, the horizontal arc connecting the right most boundary components of the surface on the right of Figure~\ref{fig:fig8}.) Therefore, the corresponding $\xi_m$ are overtwisted \cite{hkm1}. 

%We would like (in order to prove Proposition \ref{bn}) to show that $\xi_m = \xi_{(1),m}$ when $m\leq 0$. Since both contact structures are overtwisted for such $m$, it suffices to show that their 3-dimensional invariants are equal, as discussed above. This is achieved by the following lemma in combination with Lemma \ref{d3hf}.

\begin{lem}
\label{d3}
For $m\leq 0$, the 3-dimensional invariant $d_3(\xi_m)=3/2$.
\end{lem}

\begin{proof}[Proof of Lemma \ref{d3}]
Figure \ref{fig:lefsh} shows another illustration of $P$, on the left; the four topmost horizontal segments are identified with the four bottommost horizontal segments to form 1-handles. As discussed above, we can think of these curves as knots in $\#^4(S^1\times S^2)=\partial(\sharp^4(S^1\times D^3))$. The contact manifold supported by $(P,D_r^{m} \cdot \phi)$ bounds the achiral Lefschetz fibration $X$ obtained from $\sharp^4(S^1\times D^3)$ by attaching 2-handles along the curves $\gamma_2$ and $\gamma_4$ with framing $+1$, along the curves $\gamma_1$ and $\gamma_3$ with framing $-1$, and along $|m|$ parallel copies of $r$ (with respect to the blackboard framing) with framing $+1$, as indicated on the right in Figure \ref{fig:lefsh}. 

\begin{figure}[!htbp]
\labellist 
\hair 2pt 
\small
\pinlabel $r$ at 500 265
\tiny
\pinlabel $\gamma_1$ at 275 240
\pinlabel $\gamma_2$ at 210 410
\pinlabel $\gamma_3$ at 370 150
\pinlabel $\gamma_4$ at 535 150

\pinlabel $X$ at 978 520
\pinlabel $\rotatebox{180}{\reflectbox{X}}$ at 978 25
\pinlabel $Y$ at 1140 520
\pinlabel $\rotatebox{180}{\reflectbox{Y}}$ at 1140 25
\pinlabel $Z$ at 1300 520
\pinlabel $\rotatebox{180}{\reflectbox{Z}}$ at 1300 25
\pinlabel $W$ at 1460 520
\pinlabel $\rotatebox{180}{\reflectbox{W}}$ at 1460 25

\pinlabel $1/|m|$ at 1410 265
\pinlabel $-1$ at 1110 250
\pinlabel $-1$ at 1220 146
\pinlabel $+1$ at 1387 146
\pinlabel $+1$ at 1055 427

\endlabellist 
\begin{center}
\includegraphics[height = 4.3cm]{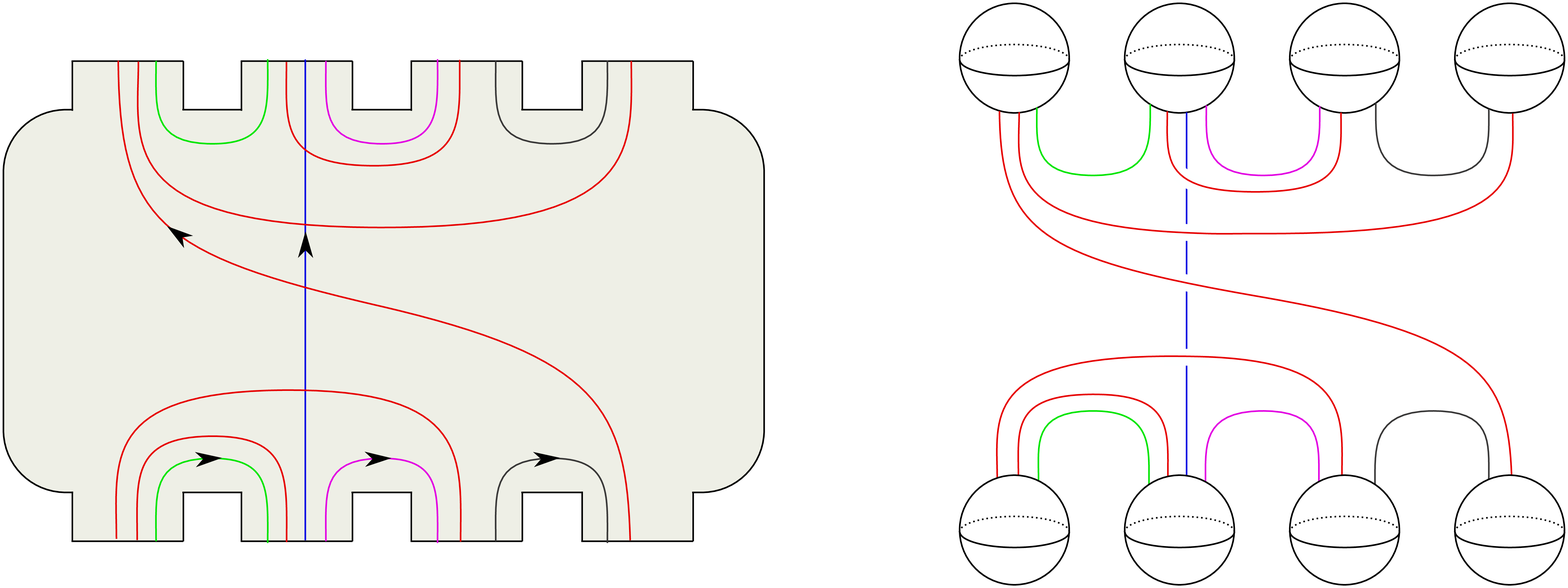}
\caption{\quad On the left, the surface $P$. On the right, a Kirby diagram for the achiral Lefschetz fibration corresponding to the open book $(P,D_r^{m} \cdot \phi)$. The label $1/|m|$ indicates that we attach 2-handles along $|m|$ parallel copies of the curve $r$ with framing $+1$.}
\label{fig:lefsh}
\end{center}
\end{figure}

Let $X$, $Y$, $Z$ and $W$ denote the 1-handles attached to $D^4$ to form $\sharp^4(S^1\times D^3)$, as shown in Figure \ref{fig:lefsh}. Furthermore, let $S_{\gamma_1}, \dots, S_{\gamma_4}$ and $S_{r_1}, \dots, S_{r_{|m|}}$ denote the cores of the 2-handles attached to the curves $\gamma_1,\dots,\gamma_4$ and the $|m|$ parallel copies $r_1,\dots, r_{|m|}$ of $r$, and let $C_{\gamma_1}, \dots, C_{\gamma_4}$ and $C_{r_1},\dots,C_{r_{|m|}}$ denote the cocores of these 2-handles. These cores form a basis for the group of 2-chains $C_2(X;\mathbb{Z})$; $X$, $Y$, $Z$ and $W$ for a basis for the 1-chains $C_1(X;\mathbb{Z})$; and the boundary map $d_2:C_2(X;\mathbb{Z})\rightarrow C_1(X;\mathbb{Z})$ sends 
\begin{eqnarray*} 
d_2(S_{\gamma_1}) &=& Y,\\
d_2(S_{\gamma_2}) &=& X-Y,\\
d_2(S_{\gamma_3}) &= &Y-Z,\\
d_2(S_{\gamma_4}) &=& Z-W,\\
d_2(S_{r_i}) &=& -Z+Y+W.
\end{eqnarray*}
The homology $H_2(X;\mathbb{Z})$ is therefore generated by $h_1,\dots,h_{|m|}$, where $$h_i = S_{r_i}+S_{\gamma_4}-S_{\gamma_1}.$$

By construction, $X$ may also be obtained from $D^4$ by attaching 2-handles along $|m|$ parallel copies of the figure eight with framing $+1$, so the intersection matrix $Q_X$ is simply the $|m|\times |m|$ identity matrix with respect to the corresponding basis. Since the curves $r_i$ are parallel, it is clear that $h_i\cdot h_j = 0$ for $i\neq j$. It follows that $h_i\cdot h_i = 1$ for $i=1,\dots,|m|$.

Recall that the class $c(X)$ is Poincar{\'e} dual to $$\sum_{i=1}^{4}\rot(\gamma_i)\cdot C_{\gamma_i}\, +\, \sum_{i=1}^{|m|}\rot(r_i)\cdot C_{r_i}.$$ Via the discussion in \cite[Section 3.1]{et3}, we calculate that $\rot(\gamma_4)=\rot(\gamma_2)=\rot(\gamma_3)=-1$ and $\rot(\gamma_1) = \rot(r_i) = 0$. Therefore, $\langle c(X),h_i\rangle = -1$ for $i=1,\dots,|m|$. So, thought of as a class in $H^2(X,\partial X;\mathbb{Z})$, $c(X)$ is Poincar\'e dual to  $$-\,h_1 \,-\, \cdots \,-\, h_{|m|}.$$ Hence, $c^2(X) = |m|$. In addition, $\chi(X) = 1+|m|$, $\sigma(X) = |m|$ and $q=2+|m|$. From the formula in (\ref{d32}), we have $$d_3(\xi_m) = \frac 14 (|m|-2(1+|m|)-3|m|) + 2 + |m| = 3/2.$$ This completes the proof of Lemma \ref{d3}. \end{proof}

\begin{figure}[!htbp]
\labellist 
\hair 2pt 
\small
\pinlabel $a$ at 225 425
\pinlabel $b$ at 350 92
\endlabellist 
\begin{center}
\includegraphics[height = 3.4cm]{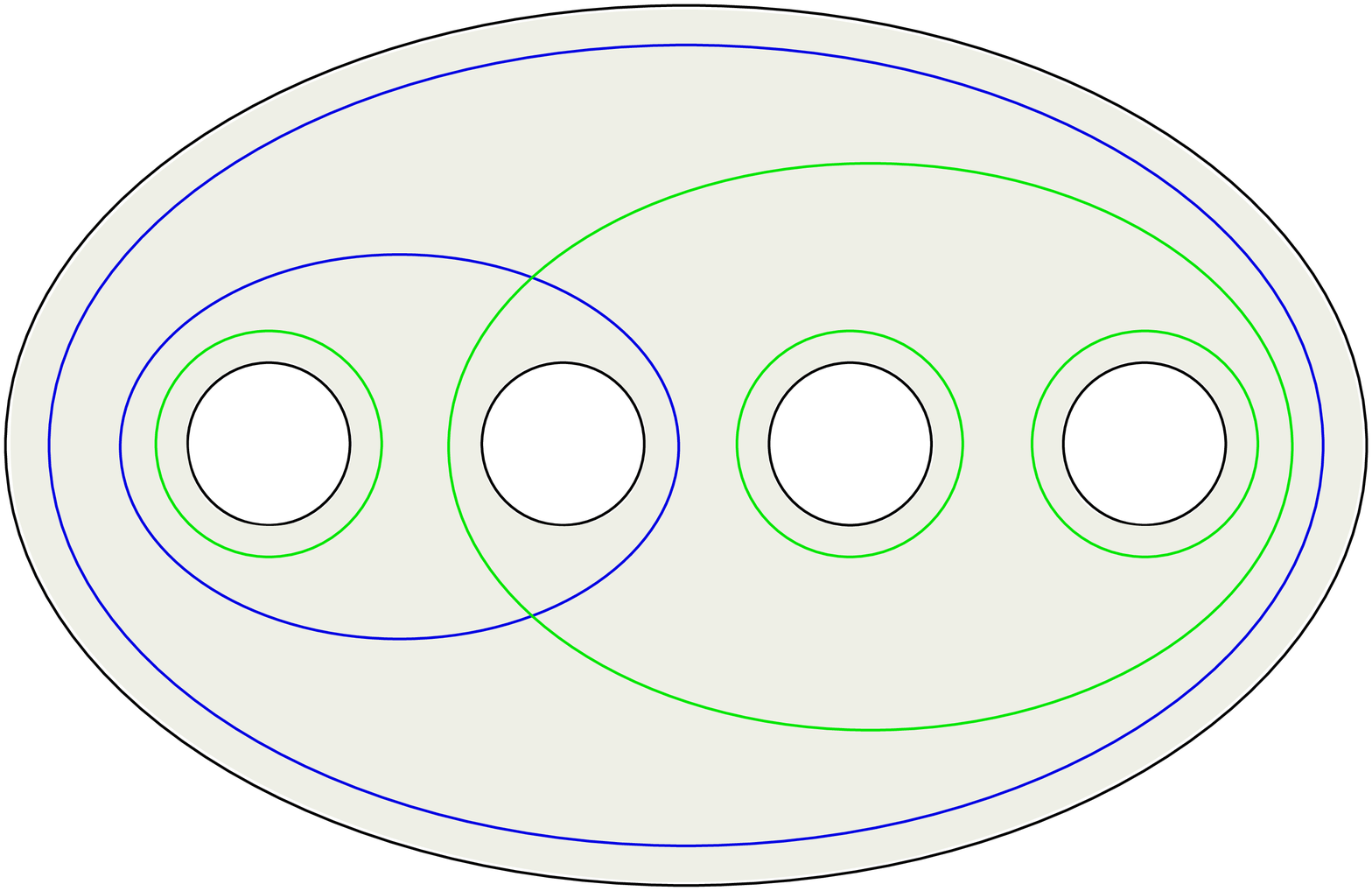}
\caption{\quad The surface $P'$.}
\label{fig:p5}
\end{center}
\end{figure}

\begin{proof}[Proof of Proposition \ref{bn}]
Recall that $$d_3(\xi\,\#\,\xi') = d_3(\xi)+d_3(\xi')+1/2$$ for any two contact structures $\xi$ and $\xi'$. Let $m<0$. Since $d_3(\xi_{(1),m}) = 1/2$, $d_3(\xi_{m}) = 3/2$, and both $\xi_{(1),m}$ and $\xi_{m}$ are overtwisted, it follows that $\xi_{(1),m}$ is isotopic to $\xi_m \,\#\, \xi',$ where $\xi'$ is the unique (overtwisted) contact structure on $S^3$ with $d_3(\xi') = -3/2$. In \cite{et3}, Ozbagci and the second author show that $\bn(\xi')\leq5$. In particular, $\xi$ is  supported by the open book $(P', D_b^{-1}D_a^{-1}\cdot\psi),$ where $a$ and $b$ are the curves on the surface $P'$ shown in Figure \ref{fig:p5} and $\psi$ is a composition of right-handed Dehn twists around the four unlabeled curves. 

Then, the planar open book $(P\,\#_b\,P', D_r^{m} \cdot \phi \cdot D_b^{-1}D_a^{-1}\cdot\psi)$ with nine binding components supports $\xi_{(1),m}\,\simeq\, \xi_m \,\#\, \xi',$ and the proof of Proposition \ref{bn} is complete. (Here, $\#_b$ denotes boundary connected sum.) \end{proof}

The table below summarizes what we know of the support genus, binding number and support norm for the contact structures $\xi_{(1),m}$. 

\begin{table}[ht]
\begin{center}
\begin{tabular}{c||c|c|c|}
 \cline{2-4} 
 
$ $&\small$\sg$& \small$\bn$&\small$\sn$ \cr\hline\hline
 \multicolumn{1}{|c||}{\small$m>0$}  & \small$1$ & \small$1$ &  \small$1$  \cr \hline 
  \multicolumn{1}{|c||}{\small$m=0$} & \small$0$ & \small$1$ & \small$-1$\cr \hline 
  \multicolumn{1}{|c||}{\small$m=-1$} & \small$0$ & \small$[3,9]$ & \small$1$ \cr \hline
   \multicolumn{1}{|c||}{\small$m<-1$} & \small$0$ & \small$[4,9]$ & \small$1$ \cr \hline
\end{tabular}
\vspace{4mm}
\caption{\quad Values of $\sg$, $\bn$ and $\sn$ for $\xi_{(1),m}$.}
\end{center}
\vspace{3mm}
\end{table}

\begin{proof}[Proof of Theorem~\ref{boundondifference}]
We first observe that $M_{\bfn,m}$ is a rational homology sphere. This can be seen by noticing that $M_{\bfn, 0}$ can be obtained as the 2--fold cover of $S^3$ branched over an alternating (non-split) link (in fact, closure of a 3--braid). Thus the determinant of the link is non-zero and hence the cardinality of the first homology of the cover is finite. Since $M_{\bfn,m}$ can be obtained from $M_{\bfn,0}$ by $1/m$ surgery on a null-homologous knot it has the same first homology. 

Let $K$ be the binding of the open book $(\Sigma, \phi_{\bfn, 0})$  in $M_{\bfn,0}.$ If we fix an overtwisted contact structure on $M_{\bfn, 0}$ we can find a Legendrian knot $L$ in the knot type $K$ with Thurston-Bennequin invariant 0 and overtwisted complement. In \cite{Celik} it was shown that there is a planar open book  $(\Sigma', \phi')$ for this overtwisted contact structure that contains $L$ on a page so that the page framing is 0.

Notice that $M_{\bfn,m}$ can be obtained from $M_{\bfn,0}$ by composing $\phi'$ with a $+1$ Dehn  twist along $m$ copies of $L$ on the page of the open book. Thus each $M_{\bfn, m}$ has an overtwisted contact structure supported by a planar open book with the same number of binding components. The number of $Spin^c$ structure on $M_{\bfn,m}$ is finite and independent of $m.$ We can get from the constructed overtwisted contact structure on $M_{\bfn, m}$ to an overtwisted contact structure realizing any $Spin^c$ structure by a bounded number of Lutz twists along generators of $H_1(M_{\bfn, m})$ all of which lie on a page of the open book. As shown in \cite{et} we may positively stabilize the open book a bounded number of times and then compose its monodromy with extra Dehn twists to achieve these Lutz twists. The number of these stabilizations depends on the number of $Spin^c$ structures on $M_{\bfn, m}$ and thus is independent of $m.$ We now have planar open books realizing overtwisted contact structures representing all $Spin^c$ structures on $M_{\bfn,m}$ with the number of binding components bounded independent of $m.$

To get an open book representing any overtwisted contact structure on $M_{\bfn,m}$ we can take these and connect sum with overtwisted contact structures on $S^3.$ In \cite{et3} it was shown that all overtwisted contact structures on $S^3$ have $\bn\leq 6.$ Thus we obtain a bound independent of $m$ on the binding number for all overtwisted contact structures on $M_{\bfn,m}$ and in particular on the $\xi_{\bfn,m}, m<0.$
\end{proof}

\begin{rem}
One can also show, in a similar manner to the proof of Theorem~\ref{boundondifference}, that the binding number of $\xi_{\bfn, m}$ is bounded by a constant depending only on the length of $\bfn.$
\end{rem}

As noted in Theorem \ref{main}, $\sg(\xi_{\bfn,m})=0$ for $m<0$ and $\sn(\xi_{\bfn,m}) =1$ for all but finitely many $m<0$; that is, these two quantities do not depend (much) on $\bfn$ or $m$. While we do know that the binding number of $\xi_{\bfn, m}$ is bounded independent of $m$ it could depend on (the length of) $\bfn.$ This suggests the following interesting question, which we leave unanswered.

\begin{quest}
Does there exist, for any positive integer $n$, a contact structure $\xi$ such that $2\sg(\xi)+\bn(\xi)-2-\sn(\xi)=n$?
\end{quest} 

Note that $0\,\leq\,2\sg(\xi_{(1),m})+\bn(\xi_{(1),m})-2-\sn(\xi_{(1),m})\,\leq\,6$ for all $m\in \mathbb{Z}$, but we currently cannot prove that this difference is larger that 1. It would be very interesting to determine if there is an $m$ such that this difference is greater than 1. In general, computing $2\sg(\xi_{\bfn,m})+\bn(\xi_{\bfn,m})-2-\sn(\xi_{\bfn,m})$ could potentially provide a positive answer to this question. 

Noticing that all our examples involve overtwisted contact structures we end with the following question.

\begin{quest}
Is there a tight contact structure $\xi$ such that $$\sn(\xi)<2\sg(\xi)+\bn(\xi)-2 \text{ ?}$$
\end{quest}

%\bibliographystyle{hplain.bst}
%\bibliography{References}

\end{document}